\documentclass[a4paper,12pt]{amsart}
\usepackage{amssymb}
\usepackage{amsmath, amsthm, amscd, amsfonts, amssymb, graphicx, color}
\usepackage{amsmath, amsthm, amscd, amsfonts, amssymb, graphicx, color}

\usepackage[cp1250]{inputenc}
\usepackage[T1]{fontenc}

\newtheorem{theorem}{Theorem}[section]
\newtheorem{lemma}[theorem]{Lemma}
\newtheorem{proposition}[theorem]{Proposition}
\newtheorem{corollary}[theorem]{Corollary}
\theoremstyle{definition}
\newtheorem{definition}[theorem]{Definition}
\newtheorem{example}[theorem]{Example}

\theoremstyle{remark}
\newtheorem{remark}[theorem]{Remark}
\numberwithin{equation}{section}

\begin{document}
	\author{Stefan Ivkovi\'{c} }

	\vspace{15pt}
	
	\title{Porosity and supercyclic operators on solid Banach function spaces}

	\maketitle
	\begin{abstract}
	In this paper, we characterize supercyclic weighted composition operators on a large class of solid Banach function spaces, in particular on Lebesgue, Orlicz and Morrey spaces. Also, we characterize supercyclic weighted composition operators on certain Segal algebras of functions and on non-unital commutative C*-algebras. Moreover, we introduce the concept of Ces\'{a}ro hyper-transitivity and we characterize Ces\'{a}ro hyper-transitive weighted composition operators on all these spaces. We illustrate our results by concrete examples and we give in addition an example of a hypercyclic weighted composition operator which is not Ces\'{a}ro hyper-transitive. Next, we introduce a class of non-porous subsets of the space of continuous functions vanishing at infinity on the real line. As an application, we consider weighted composition operator on this space and we give sufficient conditions that induce that the set of non-hypercyclic vectors for this operator is non-porous.	
	\end{abstract}
	
	\vspace{15pt}
	
	\begin{flushleft}
	\end{flushleft} 
	
	\vspace{15pt}
	
	\begin{flushleft}
	\end{flushleft}
	
	\vspace{30pt}
	
	\section{Introduction}
	
	Linear dynamics of operators have been studied in many articles during several decades; see \cite{bmbook} and \cite{gpbook} as monographs on this topic.
	Among other concepts, hypercyclicity, topological transitivity and supercyclicity have been investigated in many research works. For example,  Hilden and Wallen in \cite{Supercyclic-1} proved that any unilateral backward weighted shift is supercyclic. Afterwards, Salas in \cite{salas} characterized supercyclic bilateral weighted shift operators on $ l^p ( \mathbb{Z} ) $ in terms of a supercyclicity criterion. Supercyclicity of several kinds of operators has also been studied in for instance \cite{comparison, gupta, novo-prvi, novo-drugi, banach}. Moreover, there is a close connection between supercyclicity and semi-Fredholm theory, see \cite{Fredholm-3}. Another concept in linear dynamics of operators which is related to semi-Fredholm theory (and also to ergodic theory) is Ces\'{a}ro-hypercyclicity. This concept has been introduced and studied by Le\'{o}n-Saavedra in \cite{savedra}, see also \cite{caot-saavedra} and \cite{filomat-cesaro}. \\
	Now, the dynamics of weighted composition operators on various function spaces has attracted attention of many researchers.  In \cite{cc09,cc11},  some sufficient and necessary conditions for weighted translation operators to be hypercyclic on the Lebesgue space in context of homogeneous spaces and locally compact groups have been obtained.
	Recently, these works made a significant motivation for researchers to study the dynamics of weighted composition operators on some other function spaces such as Orlicz spaces; see \cite{cot}. Finally, chaotic and hypercyclic operators were studied on a large  class of solid Banach function spaces in \cite{solid,tabsaw}. Further, in \cite{CAOT} we characterized hypercyclic weighted composition operators on commutative non-unital C*-algebras and on certain Segal algebras of functions.
	
	We recall that supercyclicity is weaker property than hypercyclicity, i.e. every hypercyclic operator is supercyclic, whereas the converse in not true in general. The purpose of Section 3 in this paper is to extend the results from \cite{solid, CAOT, tabsaw} from the hypercyclic to the supercyclic case, thus to provide necessary and sufficient conditions for the weighted composition operators considered in \cite{solid, CAOT, tabsaw} to be supercyclic. First, motivated by the approach from \cite{solid, tabsaw}, we consider weighted composition operators on a large class of solid Banach function spaces and in Proposition \ref{proposition-solid} we characterize supercyclic such operators, with applications to Morrey spaces, see Remark \ref{applications}. The weighted composition operators on Lebesgue and Orlicz spaces are another special cases of this approach. Then we recall the construction of special Segal function algebras from \cite{tak, CAOT} and in  Proposition \ref{proposition-segal}, we characterize supercyclic weighted composition operators on these Segal algebras. Finally, in Proposition \ref{cont} we characterize supercyclic weighted composition operators on commutative non-unital C*-algebras. In this context we observe that if the weight function is non-negative, then the corresponding weighted composition operator is a completely positive operator, and the importance of completely positive operators in quantum information theory is well known. Moreover, motivated by the comparison of supercyclicity with Ces\'{a}ro-hypercyclicity for weighted translations on locally compact groups in \cite{comparison}, we introduce the concept of topological Ces\'{a}ro hyper-transitivity, and in Section 3  we characterize in addition Ces\'{a}ro hyper-transitive weighted composition operators on the above mentioned function spaces and algebras. Also, we illustrate these results by giving concrete examples of Ces\'{a}ro hyper-transitive non-hypercyclic weighted composition operators (Example \ref{semi-transitive}). By our very definition, every Ces\'{a}ro hyper-transitive operator is supercyclic, however, it turns out that the converse is not true in general. In Examples \ref{noncesaro1} and \ref{non-Cesaro-hypercyclic} we construct some supercyclic non-hypercyclic weighted composition operators that are not Ces\'{a}ro hyper-transitive. Finally, in Example \ref{hypercyclic} we construct a hypercyclic weighted composition operator on $ L^{2}(\mathbb{R}) $ which is not Ces\'{a}ro hyper-transitive, and, then, as a corollary and correction of \cite[Example 2]{pogresanprimer}, in Remark \ref{existence}, we give a new example of a hypercyclic bilateral weighted shift which is not Ces\'{a}ro hypercyclic in the sense of \cite{savedra}. At the end of Section 3, in Example \ref{compact} we show how our examples of supercyclic weighted composition operators induce examples of supercyclic left multipliers and supercyclic completely positive operators on the C*-algebra of compact operators on $ L^{2}(\mathbb{R}) $.  \\
	In  Section 4, we characterize supercyclic and Ces\'{a}ro hyper-transitive adjoints of weighted composition operators and we illustrate these results by concrete examples (Example \ref{measures}).\\
	Section 5 of this paper is devoted to the study of the dynamics of weighted composition operators in the context of the notion of porosity.  Now, $\sigma$-porous sets, as a special collection of very thin subsets of  metric spaces, were introduced and studied first time in \cite{dol}. The concepts  related to porosity have been active topics in recent decades, see the monograph \cite{zaj2}. Recently, F. Bayart   proved in \cite{bay} that the set of non-hypercyclic vectors of some classes of weighted shift operators on $\ell^2(\mathbb{Z})$ is a non-$\sigma$-porous set, and in \cite{taiwanese} this result has been extended to weighted composition operators on more general  $ L^p$-spaces. Thanks to Theorem \ref{porous} in this paper, in Corollary \ref{porosity} we obtain an extension of the results from \cite{bay,taiwanese} to the case of weighted composition operators on $ C_0( \mathbb{R}) ,$ (where $ C_0( \mathbb{R}) $ denotes as usual the space of continuous function vanishing at infinity on $ \mathbb{R} $). 
	
	\section{Preliminaries}

   To keep the paper sufficiently self-contained, we recall now the following definitions. 
  \begin{definition}
  	Let $ X$ be a separable Banach space and $B(X)$ denote the space of all linear bounded operators on $X$.  A sequence $(T_n)_{n\in \Bbb{N}}$  of operators in $B( X)$ is called {\it hypercyclic} if there is an element $x\in X$ (called \emph{hypercyclic vector}) such that the set $ \{ T_nx:\,n\in\mathbb N,  \}$ is dense in $ X$. The set of all hypercyclic vectors of a sequence $(T_n)_{n\in \Bbb{N}}$ is denoted by $HC((T_n)_{n\in \Bbb{N}})$. If $HC((T_n)_{n\in \Bbb{N}})$ is dense in $X$, the sequence $(T_n)_{n\in \Bbb{N}}$ is called \emph{densely hypercyclic}. An operator $T\in B(X)$ is called \emph{hypercyclic} if the sequence $(T^n)_{n\in \mathbb N}$ is densely
  	hypercyclic.\\
  	Further, a sequence $(T_n)_{n\in \Bbb{N}}$  of operators in $B( X)$ is called {\it supercyclic} if there is an element $x\in X$ (called \emph{supercyclic vector}) such that the set $$ \{ \lambda T_nx:\,n\in\mathbb N, \text{ } \lambda \in \mathbb{C} \setminus \lbrace 0 \rbrace  \}$$ is dense in $ X$. The set of all supercyclic vectors of a sequence $(T_n)_{n\in \Bbb{N}}$ is denoted by $SC((T_n)_{n\in \Bbb{N}})$. If $SC((T_n)_{n\in \Bbb{N}})$ is dense in $X$, the sequence $(T_n)_{n\in \Bbb{N}}$ is called \emph{densely supercyclic}. An operator $T\in B(X)$ is called \emph{supercyclic} if the sequence $(T^n)_{n\in \mathbb N}$ is densely
  	 supercyclic.
  \end{definition} 
\begin{definition}
	Let $X$ be a Banach space and $T \in B(X).$ We say that $T$ is {\it topologically transitive} on $X$ if for each pair of open non-empty subsets $O_{1}$ and $O_{2}$ of $X$ there exists some $n \in \mathbb{N}$ such that $  T^{n}(O_{1}) \cap O_{2} \neq \varnothing .$ Similarly, we say that $T$ is {\it topologically semi-transitive} on $X$ if for each pair of open non-empty subsets $O_{1}$ and $O_{2}$ of $X$ there exists some $n \in \mathbb{N}$ and some $\lambda \in \mathbb{C} \setminus \lbrace 0 \rbrace $ such that $\lambda T^{n}(O_{1}) \cap O_{2} \neq \varnothing .$ We say that $T$ is {\it topologically Ces\'{a}ro hyper-transitive} on $X$ if for each pair of open non-empty subsets $O_{1}$ and $O_{2}$ of $X$ there exists a strictly increasing sequence of natural numbers $  \lbrace n_{k} \rbrace_{k} $ such that $  n_{k}^{-1} T^{n_k}(O_{1}) \cap O_{2} \neq \varnothing $ for all $k.$
\end{definition}

In \cite[Definition 1.2]{novo-prvi}, topological semi-transitivity is actually called {\it topological transitivity for supercyclicity}. Moreover, in \cite[Proposition 1.3]{novo-prvi} it has been proved that if $X$ is a separable Banach space, then an operator $ T \in B(X) $ is topologically semi-transitive if and only if $T$ is densely supercyclic. The motivation for the concept of Ces\'{a}ro hyper-transitivity comes from \cite{savedra}.

At the end of this section, we give also the following auxiliary remark which we will use later in the proofs.  

\begin{remark} \label{remark-semi-transitive}
	
	If $T \in B(X) $ is invertible  and topologically semi-transitive,  then, given two non-empty open subsets $ O_1 $ and $ O_2 $ of $X ,$ there exists a strictly  increasing sequence $\lbrace n_{k}\rbrace_{k} \subseteq \mathbb{N}$ and sequence $\lbrace \lambda_{k} \rbrace \subseteq \mathbb{C} \setminus \lbrace 0 \rbrace $ such that for all $k \in \mathbb{N}$ we have that $$\lambda_{k} T^{n_{k}}(O_{1}) \cap O_{2} \neq \varnothing .$$ Indeed, since $T$ is topologically semi-transitive, we can find some $ n_{1} \in \mathbb{N}        $ and some  $  \lambda_{1} \in \mathbb{C} \setminus \lbrace 0 \rbrace       $  such that   $\lambda_{1} T^{n_{1}}(O_{1}) \cap O_{2} \neq \varnothing        $  . Now, since   $  \lambda_{1} T^{n_{1}}      $   is invertible, it is an open map, hence   $  \lambda_{1} T^{n_{1}}(O_{1})      $  is open. Therefore, there exists some   $ \tilde{n}_{2} \in \mathbb{N}       $  and some   $\tilde{\lambda}_{2} \in  \mathbb{C} \setminus \lbrace 0 \rbrace          $  such that   $\tilde{\lambda}_{2} T^{{\tilde{n}}_{2}}(\lambda_{1} T^{n_{1}}(O_{1})) \cap O_{2} \neq \varnothing .$   Put   $n_{2}=n_{1}+\tilde{n}_{2} $  and   $  \lambda_{2} = \lambda_{1} \tilde{\lambda}_{2} .$ Proceeding inductively, we can construct the desired sequences   $ \lbrace n_{k} \rbrace_{k} $  and   $  \lbrace \lambda_{k} \rbrace_{k} .$ 
\end{remark}

\section{Supercyclicity of weighted composition operators}

In this section, the set of all Borel measurable complex-valued functions on a topological space $X$ is denoted by $\mathcal M_0(X) .$ Also, $\chi_A$ denotes the characteristic function of a Borel set $A$. We recall the following definitions from \cite{solid}.
\begin{definition}
	Let $X$ be a topological space and $\mathcal{F}$ be a linear subspace of $\mathcal M_0(X)$. If $\mathcal F$ equipped with a given norm $\|\cdot\|_{\mathcal F}$ is a Banach space, we say that $\mathcal F$ is a \emph{Banach function space on $X$}.
\end{definition}
\begin{definition}
	Let $\mathcal{F}$ be a Banach function space on a topological space $X$, and  $\alpha:X\longrightarrow X $ be a homeomorphism. We say that $\mathcal{F}$ is \emph{$\alpha$-invariant} if for each 	$f\in \mathcal{F} $ we have 
	$f\circ \alpha^{\pm 1} \in \mathcal{F}$ and 	$\|f\circ \alpha^{\pm 1} \|_{\mathcal{F}}=\|f\|_{\mathcal{F}}$.
\end{definition}
\begin{definition}
	A Banach function space $\mathcal{F}$ on $X$ is called \emph{solid} if for each $f\in \mathcal{F}$ and $g\in\mathcal M_0(X)$, satisfying $|g|\leq |f|$, we have $g \in \mathcal{F}$
	and $\|g\|_{\mathcal{F}}\leq \|f\|_{\mathcal{F}}$.
\end{definition}
For the next results we shall also assume that the following conditions from \cite{solid} on the Banach function space $\mathcal{F}$ hold.

\begin{definition}\label{condition}
	Let $X$ be a topological space, $\mathcal F$ be a Banach function space on $X$, and $\alpha$ be an aperiodic homeomorphism of $X .$ We say that $\mathcal F$ satisfies condition $\Omega_\alpha$ if the following conditions hold:
	~\	\begin{enumerate}
		\item $\mathcal{F}$ is solid and $\alpha$-invariant;
		\item for each compact set	$E\subseteq X$ we have $\chi_{E}\in \mathcal{F}$;
		\item $\mathcal{F}_{bc}$ is dense in $\mathcal F$, where $\mathcal{F}_{bc}$ is the set of all bounded compactly supported functions in $\mathcal{F}$.
	\end{enumerate}
\end{definition}

 From now on, we shall assume that $ \alpha $ is aperiodic, that is for each compact subset $K$ of $ X$, there exists a constant $N_{K}>0$ such that for each $n\geq N_{K}$, we have $K \cap \alpha^{n}(K)=\varnothing$. For a measurable positive function $ w $ on $X ,$ we let $ w^{-1} := \frac{1}{w}  .$ If $ w  $ is a positive measurable function on $  X $ such that $ w,w^{-1}  $ are bounded, then $ \tilde{T}_{\alpha,w} $ will denote the weighted composition operator on $ \mathcal{F}  $ defined by $ \tilde{T}_{\alpha,w} (f) = w \cdot (f \circ \alpha)  $ for all $ f \in  \mathcal{F}.  $  In this case,  $ \tilde{T}_{\alpha,w} $ is invertible, and we will let $ \tilde{S}_{\alpha,w} $  denote its inverse. Under these assumptions and keeping the same notation, we provide now the following two propositions. 
\begin{proposition}\label{proposition-solid}
	The following statements are equivalent.\\
	i)  $  \tilde{ T }_{\alpha,w} $ is topologically semi-transitive on  $ \mathcal{F}.  $  \\
	ii) For each compact subset $K$ of $X,$ there exist a sequence of Borel subsets  $ \left\{ E_{k}\right\} _{k=1}^{\infty }  $  of $K$ and a strictly
	increasing sequence  $ \left\{ n_{k}\right\} _{k}\subseteq \mathbb{N}  $  such that  
	$$ \lim _{x\rightarrow \infty }\left\| \chi_{K\backslash E_{k}}\right\| _{\mathcal{F}}=0  $$  
	and 
	$$ \lim _{k\rightarrow \infty } [\left( \sup _{x\in E_{k}}  \prod ^{n_{k}-1}_{j=0}\left( w \circ \alpha ^{j}\right) ^{-1}\left( x\right) \right) \cdot 
	\left( \sup _{x\in E_{k}} \prod ^{n_{k}}_{j=1} ( w \circ \alpha ^{-j} ) ( x ) \right) ]  = 0. $$
	If $ \mathcal{F} $ is separable, then the above statements are equivalent to the fact that $  \tilde{ T }_{\alpha,w} $  is supercyclic.  Moreover, if $ \alpha $ is not aperiodic, then we only have that $ii) \Rightarrow i).$
\end{proposition}
\begin{proof}
	We prove first that  $  i) \Rightarrow  ii)  .$ Given a compact subset $ K\subseteq X  ,$  since  $ \tilde{ T }_{\alpha,w}  $  is topologically semi-transitive on  $ \mathcal{F},$ we can find a strictly increasing sequence
	$ \left\{ n_{k}\right\} _{k}\subseteq \mathbb{N}  ,$  a sequence  $ \left\{ \lambda_{k}\right\} _{k}\subseteq \mathbb{C} \setminus \lbrace 0 \rbrace  $ and a sequence  $ \left\{ f_{k}\right\} _{k} \subseteq \mathcal{F}   $ 
	such that $ \alpha^{n_{k}}\left( K \right) \cap  K  = \varnothing \text{  , } \left\| f_{k}-\chi_{k} \right\| _{\mathcal{F}}\leq \dfrac{1}{4^{k}} $  and  
	$ \parallel \lambda_{k} \tilde{ T }^{n_{k}}_{\alpha,w} (f_{k})-\chi_{K} \parallel_{ \mathcal{F} } \leq \dfrac{1}{4^{k}}  $  for all  $ k.  $ 
	This follows from Remark \ref{remark-semi-transitive}. Now, for each  $ k \in \mathbb{N}  ,$  we put
	
	$$ C_{k} = \lbrace x \in K : \vert \lambda_{k}   
	( \prod ^{n_{k}-1}_{j=0} 
	(w \circ \alpha^{j})(x) )  
	( f_{k}  \circ \alpha ^{n_{k}} ) (x) -1 \vert 
	\geq \dfrac{1}{2^{k}} \rbrace ,$$
	$$ D_{k}= \lbrace x \in K : ( \prod ^{n_{k}}_{j=1} ( w \circ \alpha^{-j} ) (x)) \vert   \lambda_{k}   \vert \text{ } \vert  f_{k} (x)  \vert \geq \frac{1}{2^{k}} \rbrace .$$ 
	
	By exactly the same arguments as in the proof of \cite[Theorem 1]{solid}, since
	$$ \parallel  \lambda_{k} ( \prod ^{n_{k}-1}_{j=0} ( w \circ \alpha^{j} ) 
	( f_{k}  \circ \alpha ^{n_{k}} ) - \chi_{k} \parallel_{ \mathcal{F} } =
	\parallel  \lambda_{k} \tilde{ T }^{n_{k}}_{\alpha,w} (f_{k})-\chi_{K} \parallel_{ \mathcal{F} } \leq \dfrac{1}{4^{k}} , $$ 
	we can get that  $ \parallel \chi_{C_{k}} \parallel_{ \mathcal{F} } \leq \dfrac{1}{2^{k}}  $ and $ \parallel  \chi_{D_{k}}   \parallel _{\mathcal{F}} < \frac{1}{2^{k}}  .$   Further, as in the proof of \cite[Theorem 1]{solid}, we let  
	
	$$
	A_{k}= \lbrace x \in K \text{ } \vert \text{ } \vert f_{k} (x)-1 \vert \geq \frac{1}{2^{k}} \rbrace ,\text{ } B_{k}= \lbrace x \in K^{c} \text{ } \vert \text{ } \vert f_{k} (x) \vert \geq \frac{1}{2^{k}} \rbrace
	.$$  
	Since, by the definition of  $ C_{k}  ,$  we have for all  $ x \in K \setminus C_{k}   $ that
	$$ 
	\left( \prod _{j=0} ^{n_{k}-1}    ( w \circ \alpha^{j} ) (x)\right) ^{-1} < \frac   {\vert \lambda_{k} \vert \text{ } \vert f_{k} \circ \alpha^{n_{k}} (x) \vert} {1-\frac{1}{2^{k}}}
	,$$ 
	by exactly the same arguments as in the proof of \cite[Theorem 1]{solid} we obtain that 
	$$
	( \prod _{j=0} ^{n_{k}-1}    ( w \circ \alpha^{j} ) (x)) )^{-1} < \frac   {\vert \lambda_{k} \vert \text{ } } {{2^{k}} -1}
	$$ 
	for all  $ x \in K \setminus (C_{k} \cup \alpha^{-n_{k}} (B_{k}))  .$ Moreover, for all $ x \in K \setminus A_k $ and each $ k \in \mathbb{N} ,$ we have that $ \vert f_{k} (x) \vert \geq 1 -\frac{1}{2^{k}} ,$ hence, for all  $ x \in K \setminus (D_{k} \cup A_{k})  ,$  we get
	$$
	\prod _{j=1} ^{n_{k}}    ( w \circ \alpha^{-j} ) (x) 
	< \frac{\frac{1}{2^{k}}}
	{\vert \lambda_{k} \vert \text{ } \vert f_{k} (x) \vert}
	<\frac{\frac{1}{2^{k}}}
	{\vert \lambda_{k} \vert (1-\frac{1}{2^{k}})} 
	= \frac{1}{\vert \lambda_{k} \vert (2^{k}-1)} 
	$$
	for each $k \in \mathbb{N} .$
	As in the proof of \cite[Theorem 1]{solid}, we put  $$ E_{k} = K \setminus (A_{k} \cup \alpha^{-n_{k}} (B_{k}) \cup C_{k} \cup D_{k})  $$ 
	and deduce that  $ \parallel \chi _{K \setminus E_{k}} \parallel_{\mathcal{F}} < \frac{4}{2^{k}}$     for all  $ k  .$ Finally, we also have that  
	$$ \left( \sup _{x\in E_{k}}  
	\prod ^{n_{k}-1}_{j=0} \left( w \circ \alpha ^{j} \right) ^{-1} \left( x\right) \right) \cdot \left( \sup _{x\in  E_{k}} 
	\prod ^{n_{k}}_{j=1} ( w \circ \alpha ^{-j} ) ( x ) \right)  
	$$ $$\leq \frac{\vert \ \lambda_{k} \vert }{2^{k}-1} 
	\cdot\frac{1}{\vert   \lambda_{k}   \vert (2^{k}-1)} = \frac{1}{(2^{k}-1)^{2}} \text{ for all } k \in \mathbb{N} .$$

	Next we prove $ii) \Rightarrow i).$ Let  $ \mathcal{O}_{1}   ,$  and  $ \mathcal{O}_{2}   $ 
	be non-empty open subsets of  $ \mathcal{F}  .$  Then we can find some  $ f \in (\mathcal{O}_{1}  \setminus \lbrace 0 \rbrace ) \cap \mathcal{F}_{b_{c}} $  and  $ g \in (\mathcal{O}_{2}  \setminus \lbrace 0 \rbrace ) \cap \mathcal{F}_{b_{c}}  $  since  $\mathcal{O}_{1}  \setminus \lbrace 0 \rbrace ) , \mathcal{O}_{2}  \setminus \lbrace 0 \rbrace )   $ are
	also open, non-empty and  $ \mathcal{F}_{b_{c}}    $  is dense in  $\mathcal{F}   .$  As in the proof of \cite[Theorem 1]{solid}, set  $ K = supp \text{ } f \cup supp \text{ } g    .$ Choose the strictly increasing sequence  $ \lbrace n_{k} \rbrace_{k} \subseteq \mathbb{N}   $ and the sequence of Borel subsets  $ \lbrace E_{k} \rbrace_{k}  $  of $K$ that satisfy the assumptions of $ii)$ with respect to $K.$ As shown in the proof of \cite[Theorem 1]{solid}, we have that  $\parallel f - f  \chi _{E_{k}} \parallel_{ \mathcal{F} } \rightarrow 0  $ and  $ \parallel g - g  \chi _{E_{k}} \parallel_{ \mathcal{F} } \rightarrow 0    $ 
	when  $ k \rightarrow \infty  ,$  so we may without loss of generality, asume that  $f  \chi _{E_{k}} \in \mathcal{O}_{1} \setminus \lbrace 0 \rbrace  $ and  $ g  \chi _{E_{k}} \in \mathcal{O}_{2} \setminus \lbrace 0 \rbrace  $  for all $k.$ Therefore,
	$  \tilde{T}_{\alpha,w}^{n_{k}}    (f  \chi _{E_{k}} ) \neq 0 $ and  $  \tilde{S}_{\alpha,w}^{n_{k}}    (g  \chi _{E_{k}} ) \neq 0  $  for all $k$ because $\tilde{T}_{\alpha,w} $ and $ \tilde{S}_{\alpha,w} $ are invertible. As in the proof of \cite[Theorem 1]{solid}, we have for each  $  k \in \mathbb{N} $  that
	
	$$ 
	\parallel 
	\tilde{T}_{\alpha,w}^{n_{k}}   ( f  \chi _{E_{k}} ) \parallel_{\mathcal{F}} \leq 
	\parallel f \parallel_{\mathcal{F}}
	\sup _{x\in E_{k}}  
	\prod _{j=1} ^{n_{k}} \left( w \circ \alpha ^{-j} \right)  \left( x\right) \text{    } (3) ,$$
	$$
	\parallel 
	\tilde{S}_{\alpha,w}^{n_{k}}   ( g  \chi _{E_{k}} ) \parallel_{\mathcal{F}} \leq 
	\parallel g \parallel_{\mathcal{F}}
	\sup _{x\in E_{k}}  
	\prod _{j=0} ^{n_{k}-1} \left( w \circ \alpha ^{j} \right)^{-1}  \left( x\right) \text{    } (4) 
	.$$
	
	Set
	$$ 
	v_{k}= f \chi _{E_{k}} + \frac
	{ \parallel   \tilde{T}_{\alpha,w}^{n_{k}}    (f  \chi _{E_{k}} )   \parallel^{\frac{1}{2}}}
	{\parallel  \tilde{S}_{\alpha,w}^{n_{k}}    (g  \chi _{E_{k}} )    \parallel^{\frac{1}{2}}}
	\tilde{S}_{\alpha,w}^{n_{k}}    (g  \chi _{E_{k}} )   
	.$$
	By combining (3) and (4) together with the assumptions in  $ii),$ it is not hard to deduce that  $ v_{k} \rightarrow f  $  and
	$$ \frac
	{ \parallel   \tilde{S}_{\alpha,w}^{n_{k}}    (g  \chi _{E_{k}} )   \parallel^{\frac{1}{2}}}
	{\parallel  \tilde{T}_{\alpha,w}^{n_{k}}    (f  \chi _{E_{k}} )    \parallel^{\frac{1}{2}}}
	\tilde{T}_{\alpha,w}^{n_{k}}    (v_{k} ) \rightarrow g   $$
	
	as  $ k \rightarrow \infty  .$ Hence, $ \tilde{T}_{\alpha,w} $ is topologically semi-transitive on  $  \mathcal{F} .$ \\
	Finally, if $  \mathcal{F} $ is separable, then the last statement of the proposition follows from \cite[Proposition 1.3]{novo-prvi}.
\end{proof}

\begin{proposition}\label{proposition-solid-Cesaro}
	The following statements are equivalent.\\
	i) $  \tilde{ T }_{\alpha,w} $ is topologically Ces\'{a}ro hyper-transitive on  $ \mathcal{F}.  $  \\
	ii) For each compact subset $K$ of $X,$ there exist a sequence of Borel subsets  $ \left\{ E_{k}\right\} _{k=1}^{\infty }  $  of $K$ and a strictly
	increasing sequence  $ \left\{ n_{k}\right\} _{k}\subseteq \mathbb{N}  $  such that  
	$$ \lim _{x\rightarrow \infty }\left\| \chi_{K\backslash E_{k}}\right\| _{\mathcal{F}}=0  $$  
	and 
	$$ \lim _{k\rightarrow \infty } \left(  \sup _{x\in E_{k}}  n_k \prod ^{n_{k}-1}_{j=0}\left( w \circ \alpha ^{j}\right) ^{-1}\left( x\right) \right) = \lim _{k\rightarrow \infty }
	\left(  \sup _{x\in E_{k}} n_{k}^{-1} 
	\prod ^{n_{k}}_{j=1} ( w \circ \alpha ^{-j} ) ( x )  \right)  = 0. $$
	Moreover, if $ \alpha $ is not aperiodic, then we only have that $ii) \Rightarrow i).$
\end{proposition}
\begin{proof}
	For the proof of the implication $  i) \Rightarrow  ii)  $ we can proceed in the similar way as in the proof of $  i) \Rightarrow  ii)  $ in Proposition \ref{proposition-solid} by letting $ \lambda_{k} = n_{k}^{-1} ,$ whereas for the proof of the opposite implication, we can proceed as in the proof of  \cite[Theorem 1]{solid} by considering $ n^{-1}  \tilde{ T }_{\alpha,w}^{n} $ and $ n  \tilde{ S }_{\alpha,w}^{n} $  instead of $  \tilde{ T }_{\alpha,w}^{n} $ and $   \tilde{ S }_{\alpha,w}^{n} ,$  respectively. 
\end{proof}	
\begin{remark}\label{applications}
For the applications of the above propositions in the case of Morrey spaces, see \cite[Example 4]{solid}.
\end{remark}

In sequel, $ \Omega$ will denote a locally compact non-compact Hausdorff space,  $ C_{0}(\Omega) $ will denote  the commutative $ C^{*}$-algebra of all continuous functions vanishing at infinity  on $ \Omega$ (equipped with the supremum norm). We let  $\mathcal{A} =C_{0}(\mathbb{R})$ and $ \tau \in C_{b}(\mathbb{R}) ,$ that is $ \tau $ is a bounded continuous function on $ \mathbb{R} .$ Put
$$\mathcal{A}_\tau:=\{f\in\mathcal{A}:\, \sum_{k=0}^\infty\|f\tau^k\|_{\infty} <\infty\}.$$
For each $f\in \mathcal{A}_\tau$ we define
$$\|f\|_\tau:=\sum_{k=0}^\infty\|f\tau^k\|_{\infty}.$$
Then, $\mathcal{A}_\tau$ is a Banach algebra \cite{tak}. Moreover, if $ \tau $ is real-valued, then $\mathcal{A}_\tau$ is a commutative Banach algebra with involution.  We will call this algebra \emph{Segal algebra corresponding to $\tau$.}
As in \cite[Section 3]{CAOT}, for any fixed $\epsilon \in (0,1),$ we shall denote by $ K_{\epsilon}^{(\tau)}  $ a (general) compact subset of $ \vert \tau \vert^{-1} ([0,\epsilon]) $. 

 If  $ w$ is continuous positive function on $ \mathbb{R}  $ such that $ w,w^{-1}  $ are bounded and $ \alpha  $ is a homeomorphism of $ \mathbb{R}  $ such that $ \tau \circ \alpha = \tau  ,$ 
then, by \cite[Lemma 3.9]{CAOT}, the operator $ \tilde{T}_{\alpha,w}  $ is an invertible bounded linear self-mapping on $ \mathcal{A}_{\tau}  .$   Under these assumptions and keeping the same notation, we give now the following two propositions.

\begin{proposition}\label{proposition-segal}
	The following statements are equivalent.\\ 
	(1) The operator $  \tilde{T}_{\alpha,w} $ is topologically  semi-transitive on $ \mathcal{A}_{\tau}  .$\\
	(2) For each positive $ \epsilon $ and every compact subset $   K_{\epsilon}^{(\tau)} \subseteq  \vert \tau \vert^{-1} ([0,\epsilon])  $ there exists a strictly increasing sequence $ \lbrace n_{k} \rbrace_{k} \subseteq \mathbb{N}  $  such that 
	$$ 
	\lim_{k \to \infty } [  \left( \sup_{t \in  K_{\epsilon}^{(\tau)}}\prod ^{n_{k}-1}_{j=0}\left( w \circ \alpha^{j-n_{k}} \right) ( t ) \right) \cdot \left( \sup_{t \in  K_{\epsilon}^{(\tau)}}
	\prod ^{n_{k}-1}_{j=0}\left(  w \circ \alpha^{j}   \right)^{-1} ( t ) \right) ] = 0
	.$$
	If $ \mathcal{A}_{\tau} $ is separable, then the above statements are equivalent to the fact that $  \tilde{T}_{\alpha,w} $ is supercyclic on $ \mathcal{A}_{\tau} .$ Moreover, if $ \alpha $ is not aperiodic, then we only have that $  (2) \Rightarrow (1).  $
\end{proposition}
\begin{proof}
	Assume that (1) holds. Let $ \epsilon_{1}, \epsilon_{2} \in  (0,1) $ with $ \epsilon_{2} < \epsilon_{1} $ and $K_{\epsilon_{2}}^{(\tau)}$ be a compact subset of $ \vert \tau \vert^{-1} ([0,\epsilon_{2}]) .$ By \cite[Lemma 3.2]{CAOT} we can choose a function  $  \mu_{ K_{\epsilon_{2},\epsilon_{1}}^{(\tau)} }  \in \mathcal{A}_{\tau}$ satisfying that  
	$ \mu_{K_{\epsilon_{2},\epsilon_{1}}^{(\tau)} } =1  $ on  $  K_{\epsilon_{2}}^{(\tau)} $ and that
	$ supp \text{ } \mu_{K_{\epsilon_{2},\epsilon_{1}}^{(\tau)} }  $  is compact subset of  $ \vert \tau \vert^{-1} ([0, \epsilon_{1}])   .$  Hence, we shall denote  $ supp \text{ } \mu_{ K_{\epsilon_{2},\epsilon_{1}}^{(\tau)} } 
	$ by $K_{\epsilon_{1}}^{(\tau)} . $ Since  $ \alpha $  is aperiodic,
	we can find some  $ n_{1} \in \mathbb{N}  $  such that  $ \alpha^{n_{1}} (K_{\epsilon_{1}}^{(\tau)} ) \cap K_{\epsilon_{1}}^{(\tau)} = \varnothing  .$  It follows that  $ \alpha^{n_{1}} (K_{\epsilon_{2}}^{(\tau)} ) \cap K_{\epsilon_{1}}^{(\tau)} = \varnothing  $  since  $K_{\epsilon_{2}}^{(\tau)} \subseteq K_{\epsilon_{1}}^{(\tau)}   .$  
	Now, since  $ \tilde{T}_{\alpha,w}  $ is topologically semi-transitive, we can find some  $f \in \mathcal{A}_{\tau}   $ and some $ \lambda_{1} \in \mathbb{C} \setminus \lbrace 0 \rbrace $ such that  $\parallel f- \mu_{ K_{\epsilon_{2},\epsilon_{1}}^{(\tau)} } \parallel_{\tau} < \frac{1}{2}   $ and  $ \parallel \lambda_{1} \tilde{T}_{\alpha,w}^{n_{1}}    (f)- \mu_{ K_{\epsilon_{2},\epsilon_{1}}^{(\tau)} } \parallel_{\tau} < \frac{1}{4}    .$  Since  $ \parallel \cdot \parallel_{\tau} \geq \parallel \cdot \parallel_{\infty}  ,$  we
	obtain that  $ \parallel f- \mu_{ K_{\epsilon_{2},\epsilon_{1}}^{(\tau)} } \parallel_{\infty} < \frac{1}{4}   $  and  
	
	$$ \parallel ( \lambda_{1} \prod ^{n_{1}-1}_{j=0}( w \circ \alpha^{j} ) ) (f \circ \alpha^{n_{1}}) -  \mu_{ K_{\epsilon_{2},\epsilon_{1}}^{(\tau)} } \parallel_{\infty} < \frac{1}{4}  .$$ 
	By exactly the same arguments as in the proof of \cite[Theorem 2.7]{CAOT} we can deduce that  
	$$ \sup_{t \in  K_{\epsilon_{2}}^{(\tau)}} \vert \lambda_{1} \vert \prod ^{n_{1}-1}_{j=0}( w \circ \alpha^{j-n_{1}} ) (t)  < \frac{1}{2} $$
	and
	$$  \sup_{t \in  K_{\epsilon_{2}}^{(\tau)}}  \frac{1}{ \vert\lambda_{1} \vert } \prod ^{n_{1}-1}_{j=0}( w \circ \alpha^{j} )^{-1} (t)  < \frac{2}{3}   .$$  
	Therefore, 
	$$ \left( \sup_{t \in  K_{\epsilon_{2}}^{(\tau)}}   \prod ^{n_{1}-1}_{j=0}( w \circ \alpha^{j-n_{1}} ) (t)   \right) \cdot \left(\sup_{t \in  K_{\epsilon_{2}}^{(\tau)}}  \prod ^{n_{1}-1}_{j=0}( w \circ \alpha^{j} )^{-1} (t) \right) < \frac{1}{3} .$$ 
	
	Next, we can find some  $n_{2} > n_{1}   ,$  some  $ \lambda_{2} \in \mathbb{C} \setminus \lbrace 0 \rbrace   $ and some  $ g \in \mathcal{A}_{\tau}   $ such that  $\alpha^{n_{2}}   ( K_{\epsilon_{1}}^{(\tau)} ) \cap  K_{\epsilon_{1}}^{(\tau)}  = \varnothing ,$ 
	$\parallel g - \mu_{ K_{\epsilon_{2},\epsilon_{1}}^{(\tau)} } \parallel_{\tau} < \frac{1}{2^{2}} $ 
	and  
	$ \parallel \lambda_{2} \tilde{T}_{\alpha,W}^{n_{2}} (g) - \mu_{ K_{\epsilon_{2},\epsilon_{1}}^{(\tau)} } \parallel_{\tau} <\frac{1}{4^{2}} .$ Then, as above, we can conclude that  
	
	$$ \left( \sup_{t \in  K_{\epsilon_{2}}^{(\tau)}}   \prod ^{n_{2}-1}_{j=0}( w \circ \alpha^{j-n_{1}} ) (t)   \right)\cdot \left(\sup_{t \in  K_{\epsilon_{2}}^{(\tau)}}  \prod ^{n_{2}-1}_{j=0}( w \circ \alpha^{j} )^{-1} (t) \right) < \frac{1}{4^{2}-1}   .$$   
	
	Proceeding inductively, we can construct a strictly increasing sequence  $ \lbrace n_{k} \rbrace_{k}  \subseteq \mathbb{N}  $ satisfying the assumptions in (2) with respect to  $ K_{\epsilon_{2}}^{(\tau)}, $ so the implication  $ (1) \Rightarrow (2)  $  follows.
	
	Suppose now that $(2)$ holds. Given two non-empty open subsets  $ \mathcal{O}_{1} $  and  $  \mathcal{O}_{2}  $  of 	$ \mathcal{A}_{\tau}  ,$  we can find some  $ f \in \mathcal{O}_{1}  \setminus \lbrace 0 \rbrace  $ and  $ g \in \mathcal{O}_{2}  \setminus \lbrace 0 \rbrace    .$ Then $ \tilde{T}_{\alpha,w} (f) \neq 0 $ and $ \tilde{S}_{\alpha,w} (g) \neq 0$ because $ \tilde{T}_{\alpha,w} $ and $ \tilde{S}_{\alpha,w} $ are invertible. By \cite[Corollary 3.4]{CAOT} we may assume that  $supp \text{ } f   $ and  $ supp \text{ } g   $ are compact and contained in  $\vert \tau \vert^{-1}([0,\epsilon])   $ for some  $ \epsilon \in (0,1) .$  
	Then  $  supp \text{ } f \cup supp \text{ } g  $ is also a compact subset of  $\vert \tau \vert^{-1}([0,\epsilon])   .$  Put  $ K_{\epsilon}^{(\tau)} = supp \text{ } f \cup supp \text{ } g $ and choose the strictly increasing sequence  $ \lbrace n_{k} \rbrace_{k}   $ satisfying the assumptions of $(2)$ with respect to  $ K_{\epsilon}^{(\tau)}   .$ \\
	For each  $ k \in \mathbb{N}  ,$  set
	$$ v_{k} = f + \sqrt{\dfrac{\parallel \tilde{T}_{\alpha,w}^{n_{k}} (f) \parallel_{\tau}}{\parallel \tilde{S}_{\alpha,w}^{n_{k}} (g) \parallel_{\tau}}}  \tilde{S}_{\alpha,w}^{n_{k}} (g) .$$ By the proof of \cite[Theorem 3.10]{CAOT} we have
	$$  \parallel \tilde{T}_{\alpha,w}^{n_{k}} (f) \parallel_{\tau} \leq \sup_{t \in K_{\epsilon}^{(\tau)} }  \prod ^{n_{k}-1}_{j=0}( w \circ \alpha^{j-n_{k}} ) (t)  \parallel f \parallel_{\tau }$$ 
	and
	$$  \parallel \tilde{S}_{\alpha,w}^{n_{k}} (g) \parallel_{\tau} \leq \sup_{t \in K_{\epsilon}^{(\tau)} }  \prod ^{n_{k}-1}_{j=0}( w \circ \alpha^{j} )^{-1} (t)  \parallel g \parallel_{\tau }$$
	for each $ k \in \mathbb{N} .$ Therefore, by the assumptions in (2), it follows that $v_{k} \rightarrow f  $ and
	$$ \sqrt{\dfrac{\parallel \tilde{S}_{\alpha,w}^{n_{k}} (g) \parallel_{\tau}}{\parallel \tilde{T}_{\alpha,w}^{n_{k}} (f) \parallel_{\tau}}} \tilde{T}_{\alpha,w}^{n_{k}} (v_{k}) \rightarrow g  $$
	as $k \rightarrow \infty  $ in $ \mathcal{A}_{\tau} .$	Hence, we conclude that  $ \tilde{T}_{\alpha,w}   $ is topologically semi-transitive on  $  \mathcal{A}_{\tau} .$\\
	Finally, if $  \mathcal{A}_{\tau} $ is separable, then the last statement of the proposition follows from \cite[Proposition 1.3]{novo-prvi}. 
\end{proof}
\begin{proposition}\label{proposition-segal-cesaro}
	The following statements are equivalent.\\ 
	(1) The operator $  \tilde{T}_{\alpha,w} $ is topologically  Ces\'{a}ro hyper-transitive on $ \mathcal{A}_{\tau}  .$\\
	(2) For each positive $ \epsilon $ and every compact subset $   K_{\epsilon}^{(\tau)} \subseteq  \vert \tau \vert^{-1} ([0,\epsilon])  $ there exists a strictly increasing sequence $ \lbrace n_{k} \rbrace_{k} \subseteq \mathbb{N}  $  such that 
	$$ 
	\lim_{k \to \infty }   \left( \sup_{t \in  K_{\epsilon}^{(\tau)}} n_{k}^{-1} \prod ^{n_{k}-1}_{j=0}\left( w \circ \alpha^{j-n_{k}} \right) ( t ) \right) = \lim_{k \to \infty }  \left( \sup_{t \in  K_{\epsilon}^{(\tau)}} n_k
	\prod ^{n_{k}-1}_{j=0}\left(  w \circ \alpha^{j}   \right)^{-1} ( t ) \right)  = 0
	.$$
	If  $ \alpha $ is not aperiodic, then we only have that $  (2) \Rightarrow (1).  $
\end{proposition}
\begin{proof}
	For the proof of the implication $  (1) \Rightarrow  (2)  $ we can proceed in the similar way as in the proof of $  i) \Rightarrow  ii)  $ in Proposition \ref{proposition-segal} by letting $ \lambda_{k} = n_{k}^{-1} ,$ whereas for the proof of the opposite implication, we can proceed as in the proof of  \cite[Theorem 3.10]{CAOT} by considering $ n^{-1}  \tilde{ T }_{\alpha,w}^{n} $ and $ n  \tilde{ S }_{\alpha,w}^{n} $  instead of $  \tilde{ T }_{\alpha,w}^{n} $ and $   \tilde{ S }_{\alpha,w}^{n} ,$  respectively. 
\end{proof}	
 Recall that if $ \mathcal{C} $ is a non-unital commutative C*-algebra, then by Gelfand-Naimark theorem, $ \mathcal{C} $ is isometrically *-isomorphic to $ C_{0}(\Omega)$ for some locally compact, non-compact Hausdorff space $ \Omega .$ In this case, if $w$ is a positive continuous bounded function on $ \Omega ,$ then $ \tilde{T}_{\alpha,w}   $ is a completely positive operator on  $ C_{0}(\Omega) .$
Under the additional assumption that $ \alpha $ is aperiodic and that $ w^{-1} $ is bounded, similarly as in Proposition \ref{proposition-segal} and Proposition \ref{proposition-segal-cesaro}, we can prove the following two propositions regarding the dynamics of the completely positive operator $ \tilde{T}_{\alpha,w}   $ on $ C_{0}(\Omega) .$  

\begin{proposition}\label{cont}
	The following statements are equivalent. \\
	(1) The completely positive operator $  \tilde{T}_{\alpha,w} $ is topologically  semi-transitive on $ C_{0}(`\Omega) .$\\
	(2) For every compact subset $ K $ of $ X $ there exists a strictly increasing sequence $ \lbrace n_{k} \rbrace_{k} \subseteq \mathbb{N} $ such that 
	
	$$ \lim_{k \to \infty } 
	[	\left( \sup_{t \in  K}\prod ^{n_{k}-1}_{j=0}\left( w \circ \alpha^{j-n_{k}} \right) ( t ) \right) \cdot      
	\left( \sup_{t \in  K}
	\prod ^{n_{k}-1}_{j=0}\left(  w \circ \alpha^{j}   \right)^{-1} ( t ) \right) ] = 0 .$$ 
	
	If $ C_{0}(\Omega) $ is separable, then the above statements are equivalent to the fact that $   \tilde{T}_{\alpha,w} $ is supercyclic on $ C_{0}(\Omega).$ Moreover, if $ \alpha $ is not aperiodic, then we only have that $  (2) \Rightarrow (1).  $	
\end{proposition}

\begin{proposition}\label{cont-Cesaro}
	The following statements are equivalent. \\
	(1) The completely positive operator $  \tilde{T}_{\alpha,w} $ is topologically  Ces\'{a}ro hyper-transitive on $ C_{0}(\Omega) .$\\
	(2) For every compact subset $ K $ of $ X $ there exists a strictly increasing sequence $ \lbrace n_{k} \rbrace_{k} \subseteq \mathbb{N} $ such that 
	
	$$ \lim_{k \to \infty } 
		\left( \sup_{t \in  K} n_{k}^{-1} \prod ^{n_{k}-1}_{j=0}\left( w \circ \alpha^{j-n_{k}} \right) ( t ) \right) = \lim_{k \to \infty }      
	\left( \sup_{t \in  K} n_k
	\prod ^{n_{k}-1}_{j=0}\left(  w \circ \alpha^{j}   \right)^{-1} ( t ) \right)  = 0 .$$ 
	
	If $ \alpha $ is not aperiodic, then we only have that $  (2) \Rightarrow (1).  $	
\end{proposition}
 Motivated by \cite[Example 3.6]{savedra}, we give now an example of Ces\'{a}ro hyper-transitive weighted composition operators.
\begin{example}\label{semi-transitive}
	Let $ \Omega = \mathbb{R}, \text{ } \alpha : \mathbb{R} \rightarrow \mathbb{R}   $  be given by  $\alpha (t) = t-1  $ for all  $ t \in \mathbb{R} $ and $w $  be a continuous bounded positive function on  $ \mathbb{R} .$  
	If there exist some  $ M, \delta , K_{1} , K_{2} > 0 $  such that  $ 1 < M-\delta \leq w(t) \leq M $  for all  $ t \leq -K_{1}  $  and  $ w(t) = 1  $  for all  $t \geq K_{2}  ,$  then the conditions of Proposition \ref{proposition-solid-Cesaro} and Proposition \ref{cont-Cesaro} are satisfied.
	 
\end{example}

Notice that by definition every topologically Ces\'{a}ro hyper-transitive operator is topologically semi-transitive (supercyclic). Observe also that a bounded invertible linear operator is topologically (semi-)transitive if and only if its inverse is topologically (semi-)transitive. Now we give an example of a bounded invertible supercyclic weighted composition operator which is not Ces\'{a}ro hyper-transitive, but whose inverse is Ces\'{a}ro hyper-transitive.
\begin{example}\label{noncesaro1}
	Let again $ \Omega = \mathbb{R}$ and $  \alpha : \mathbb{R} \rightarrow \mathbb{R}   $  be given by  $\alpha (t) = t +1  $ for all  $ t \in \mathbb{R} .$ If $w $  is a continuous bounded positive function on  $ \mathbb{R} $  
	such that there exist some constants $ M, \delta , K_{1} , K_{2} > 0 $
	with $  \frac{1}{M} \leq w(t) \leq \frac{1}{M- \delta}  $ for all  $t \leq -K_{1} $ and $1 = w(t) $ for all $ t\geq K_{2},$ then the conditions of Proposition \ref{proposition-solid} and Proposition \ref{cont}  are satisfied, however, by Proposition \ref{proposition-solid-Cesaro} and Proposition \ref{cont-Cesaro} it follows that  $ \tilde{T}_{\alpha, w} $ is not topologically Ces\'{a}ro hyper-transitive neither on $ L^{2}(\mathbb{R}) $ nor on $ C_{0}(\mathbb{R}) .$  However, $ \tilde{S}_{\alpha, w} $ will be topologically Ces\'{a}ro hyper-transitive both on $ L^{2}(\mathbb{R}) $ and on $ C_{0}(\mathbb{R}) $, which can be proved by considering $ n^{-1}  \tilde{ S }_{\alpha,w}^{n} $ and $ n  \tilde{ T }_{\alpha,w}^{n} $ instead of $ n^{-1}  \tilde{ T }_{\alpha,w}^{n} $ and $ n  \tilde{ S }_{\alpha,w}^{n} ,$ respectively, in the proof of the implication $ ii) \Rightarrow i) $ in Proposition \ref{proposition-solid-Cesaro}. 
\end{example}
In the next example, we shall construct  a bounded invertible supercyclic weighted composition operator which is not Ces\'{a}ro hyper-transitive and whose inverse is neither Ces\'{a}ro hyper-transitive.

\begin{example}\label{non-Cesaro-hypercyclic}
	Let $ M, \delta $ be positive constants such that $ M \geq 2 + 2 \delta $  and $ \delta \geq 1 .$  Put $ \alpha $ to  be the function on $ \mathbb{R} $ given by $\alpha (t) = t-1  $ for all  $ t \in \mathbb{R} .$ Set \\
	$w(t) =
	\begin{cases} 
		M \text{ for }  t \leq -1, \\
		M +	\frac{t+1}{2} ( 1 + \delta - M)  \text{ for } t \in [-1, 1], \\
		1+ \delta \text{ for } t \geq 1. \\
		
	\end{cases}
	.$\\

Noticing that $ \tilde{S}_{\alpha, w} = (w \circ \alpha ^{-1})^{-1} \cdot ( f \circ \alpha ^{-1}) $ for all $ f \in C_{0}(\mathbb{R})$ and all $ f \in L^{2}(\mathbb{R}) $ (depending on whether we consider $ \tilde{T}_{\alpha, w} $ as an operator on $ C_{0}(\mathbb{R}) $ or on $ L^{2}(\mathbb{R}) $), by the above propositions applied both to $ \tilde{T}_{\alpha, w} $ and to $ \tilde{S}_{\alpha, w} $ it is not hard to deduce by some calculations that both $ \tilde{T}_{\alpha, w} $ and $ \tilde{S}_{\alpha, w} $ are supercyclic, but not Ces\'{a}ro hyper-transitive operators on  $ L^{2}(\mathbb{R}) $ and on $ C_{0}(\mathbb{R}) $. 
\end{example}

\begin{remark}Similar conclusions hold if we consider $  \mathcal{A}_{\tau} $ instead of  $ L^{2}(\mathbb{R}) $ and $ C_{0}(\mathbb{R}) $ in the above examples, provided that in this case $ \tau $ satisfies that $ \tau (t+1) = \tau (t) $ for all $ t \in \mathbb{R} $ (because then $ \tau \circ \alpha = \tau  ,$ so $ \tilde{T}_{\alpha,w} $ is well defined operator  on $  \mathcal{A}_{\tau} $ by \cite[Lemma 3.9]{CAOT}).
	\end{remark}
 It follows from \cite[Theorem 1]{solid} and \cite[Corollary 2.8]{CAOT} that the weighted composition operators constructed in Example \ref{semi-transitive}, Example \ref{noncesaro1} and Example \ref{non-Cesaro-hypercyclic} are not  hypercyclic neither on $ L^{2}(\mathbb{R}) $ nor on $ C_{0}(\mathbb{R}) .$ Now we will provide an example of a hypercyclic weighted composition operator on $ L^{2}(\mathbb{R}) $ which is not Ces\'{a}ro hyper-transitive.
\begin{example}\label{hypercyclic}
Put again $ \alpha $ to be the function on $ \mathbb{R} $ given by $\alpha (t) = t-1  $ for all  $ t \in \mathbb{R} .$	Let $ w $ be the function on $ \mathbb{R} $ given by  $ w(t) = \frac{1}{2} $ for all $ t \geq 0,$ and for each $ m \in \mathbb{N} $ and $ t \in [-m, -m+1)  $  be given by $ w(t) = \frac{m+1}{m} $. If we put $ \mathcal{F} = L^{2}(\mathbb{R})
 ,$ then it is not hard to check by some calculations that the conditions of \cite[Theorem 1]{solid} are satisfied, whereas the conditions of  Proposition \ref{proposition-solid-Cesaro} are not satisfied, hence the corresponding weighted composition operator $ \tilde{T}_{\alpha,w} $ will be a hypercyclic, but not Ces\'{a}ro hyper-transitive operator on $ L^{2}(\mathbb{R}) .$ 
\end{example} 
\begin{remark}\label{existence}
	In \cite[Example 2]{pogresanprimer}, it is claimed that  the bilateral backward weighted shift on $ l^{2}( \mathbb{Z})  $ with the weight sequence $ \lbrace w_j \rbrace$ given by $ w_j = 2 $ for all $ j < 0 $ and $ w_j = \frac{1}{2} $  for all $ j \geq 0 $ is hypercyclic, but not Ces\'{a}ro hypercyclic. For hypercyclicity of this operator, the authors in \cite{pogresanprimer} refer to \cite[Theorem 4.1]{feldman}, however, \cite[Theorem 4.1]{feldman} deals with forward shifts (and not backward shifts). From \cite[Theorem 2.1]{sa95} it is straightforward to check that the operator from \cite[Example 2]{pogresanprimer} is not hypercyclic. However, if we let, as in Example \ref{hypercyclic}, $ T $ be the bilateral weighted forward shift operator with weights $ w_{-j} = \frac{j+1}{j} $ for all $ j \in \mathbb{N}$ and $ w_j = \frac{1}{2}$ for all $j \geq 0 ,$ then $T$ would be hypercyclic, but not Ces\'{a}ro hypercyclic.
\end{remark}
For examples of hypercyclic Ces\'{a}ro hyper-transitive weighted composition operators, we refer to for instance \cite[Example 2.10]{CAOT}.\\
At the end of this section we wish to illustrate how the above examples of supercyclic weighted composition operators induce examples of supercyclic operators on the C*-algebra $ K(L^{2} (\mathbb{R})) $ of compact operators on $ L^{2} (\mathbb{R}) .$  

\begin{example}\label{compact}
	For  an invertible operator $ W \in B(L^{2} (\mathbb{R})) $ and a unitary operator $ U $ on  $L^{2} (\mathbb{R}) ,$ we will denote by $ C_W$ the completely positive operator on $ K(L^{2} (\mathbb{R})) $ given by $ C_W (F) = WFW^{*} ,$ and by $ T_{U,W} $  the wedge operator on $ K(L^{2} (\mathbb{R})) $ given by $ T_{U,W} (F) = WFU$  for all $ F \in K(L^{2} (\mathbb{R})) $ (notice that $ T_{I,W} $ is in fact the left multiplier by $ W$). Let  now $W$ be the operator on $L^{2} (\mathbb{R})$ defined by  $ W(f)=w \cdot (f \circ \alpha ) $  for all  $ f \in L^{2} (\mathbb{R}) ,$ where $ \alpha $ is a homeomorphism of $ \mathbb{R} $ and $ w $ is a measurable, bounded positive function on $\mathbb{R} $ satisfying that $ w^{-1} $ is also bounded. Then $W$ is a bounded invertible linear operator on $L^{2} (\mathbb{R}).$ If  $ m \in \mathbb{N} $  and  $f \in L^{2} (\mathbb{R})  $  with  $ supp \text{ } f \subseteq [-m,m] ,$  then, by some calculations, it is not hard to see that for all  $ n \in \mathbb{N}  $ we have that
	
	$$  \int \left| W^{-n}\left( f\right) \right| ^{2}d\mu \leq
	\sup _{t \in \left[ -m,m \right] } 
	( \prod^{n-1}_{j=0}
	( w \circ \alpha^{j} )^{-1} ( t ) )  ^{2}
	\left\| f\right\|_{2}^{2} $$
	and
	$$ 
	\int \left| W^{n}\left( f\right) \right| ^{2}d\mu \leq
	\sup _{t \in \left[ -m,m \right] } 
	( \prod^{n-1}_{j=0}
	( w \circ \alpha^{j-n} ) ( t ) )  ^{2}
	\left\| f\right\|_{2}^{2} 
	.$$

	Hence, if $ w$ and $ \alpha $ satisfy that 
	$$
	\lim_{n \to \infty } [ \left( \sup _{t \in \left[ -m,m \right] }	 \prod^{n-1}_{j=0}  ( w \circ \alpha^{j-n} ) ( t ) \right) \cdot
	\left(\sup _{t \in \left[ -m,m \right] } 
	\prod^{n-1}_{j=0}
	( w \circ \alpha^{j} )^{-1} ( t ) \right) ] 
	=0, 
	$$
	for every  $ m \in \mathbb{N} ,$  then we obtain that $$  \lim_{n \to \infty } (\parallel W^{n} (f) \parallel_{2} \cdot \parallel W^{-n} (g) \parallel_{2}) = 0  $$ for every $ f,g \in  L_{c} (\mathbb{R}) ,$ where $ L_{c} (\mathbb{R}) $ denotes the space of all compactly supported functions in $ L^{2} (\mathbb{R}) .$  Since $ L_{c} (\mathbb{R}) $ is dense in $ L^{2} (\mathbb{R}) ,$ by \cite[Theorem 3.2]{gupta} we obtain that the operators $C_{W}  $ and  $ T_{U,W} $ are supercyclic on $ K(L^{2} (\mathbb{R})) $ for every unitary operator $U $ on $L^{2} (\mathbb{R}) .$ In particular,  if we let  $\alpha $ and $w $ be as  in the above examples, then these conditions are satisfied. 
\end{example}

\section{Supercyclic adjoints of weighted composition operators}

When $ \tilde{ T }_{\alpha,w}  $ is considered as an operator on $ C_{0}(\Omega) ,$ then the adjoint $ \tilde{ T }_{\alpha,w}^{*}$ is a bounded linear operator on $M(\Omega) ,$ where $M(\Omega) $ stands for the Banach space of all Radon measures on $\Omega $ equipped with the total variation norm. It is straightforward to check that $$ \tilde{ T }_{\alpha,w}^{*}(\mu)(E) = \int_{E} w \circ \alpha^{-1} \,d\mu\ \circ \alpha^{-1} $$ 
for every $\mu \in M(\Omega) $ and every measurable subset $E$ od $\Omega .$ Here $\mu\ \circ \alpha^{-1} (E) = \mu ( \alpha^{-1} (E)) $ for every  $\mu \in M(\Omega) $ and every measurable subset $E$ od $\Omega .$ Then it is not hard to check that 
$$\tilde{ T }_{\alpha,w}^{*n}(\mu)(E) = \int_{E}  \prod_{j=0}^{n-1} w \circ \alpha^{j-n} \,d\mu\ \circ \alpha^{-n} $$ 
and 
$$ \tilde{S}_{\alpha,w}^{*n}(\mu)(E) = \int_{E} \prod_{j=1}^{n} (w \circ \alpha^{n-j})^{-1} \,d\mu\ \circ \alpha^{n}.$$ 

In the rest of this section, for every Radon measure $ \mu $ on $ \Omega ,$ we let as usual $ \vert \mu \vert $ denote the total variation of $ \mu .$ Also, we assume as before that $ \alpha $ is an aperiodic homeomorphism of $ \Omega .$ Under these assumptions and keeping this notation, we obtain the following two propositions. 
\begin{proposition}\label{radon}
	The following statements are equivalent.\\
	i)  $\tilde{ T }_{\alpha,w}^{*}  $ is topologically semi-transitive on  $ M(\Omega) .$ \\
	ii) For every compact subset  $ K $ of $ \Omega $ and any two measures  $\mu, v  $ in  $ M(\Omega)  $  with  $ \vert \mu \vert  (K^{c}) = \vert v \vert  (K^{c}) =0 $ there exist a strictly increasing sequence  $ \lbrace n_{k} \rbrace_{k}  \subseteq \mathbb{N} $  and sequences  $ \lbrace A_{k} \rbrace_{k} , \lbrace B_{k} \rbrace_{k} $ of Borel subsets of $K$ such that
	$ \alpha^{n_{k}}(K) \cap K = \varnothing $  for all  $ k \in \mathbb{N}  $  and
	$$\lim_{k \to \infty } \vert \mu \vert (A_{k}) = \lim_{k \to \infty } \vert v \vert (B_{k})=0,  $$
	
	$$ \lim_{k \to \infty } 
	[	\left( \sup_{t \in A^{c}_{k} \cap K}\prod ^{n_{k}-1}_{j=0}\left( w \circ \alpha^{j}\right) ( t ) \right) \cdot      
	\left( \sup_{t \in B^{c}_{k} \cap K}
	\prod ^{n_{k}}_{j=1}\left( w \circ \alpha^{-j}\right)^{-1} ( t ) \right) ] =0 .$$ 
	If $ M(\Omega) $ is separable, then $\tilde{ T }_{\alpha,w}^{*}  $ is supercyclic in this case. Moreover, if $ \alpha $  is not aperiodic, then we only have $ ii) \Rightarrow i) .$
\end{proposition}
\begin{proof}
	We will prove $i)=> ii)$ first. Suppose that  $ \tilde{T}_{\alpha,w}^{*} $  is topologically semi-transitive on  $ M(\Omega) .$  For a given compact subset $K$ of  $\Omega  ,$  choose some
	$ \mu , v \in  M(\Omega)  $  with  $ \vert \mu \vert (K^{c}) = \vert v \vert (K^{c})= 0 .$  Similarly as in the proof of \cite[Proposition 3.1]{Grada}, by Remark \ref{remark-semi-transitive} we can find a 	sequence  $ \lbrace \eta^{(k)} \rbrace_{k} \subseteq  M(\Omega) ,$ a strictly increasing sequence  $ \lbrace n_{k} \rbrace_{k} \in \mathbb{N} $ and a sequence $\lbrace \lambda_{k} \rbrace_{k} \subseteq \mathbb{C} \setminus \lbrace 0 \rbrace$ such that for all $k$ we have that  
	$$ \alpha ^{n_k} (K) \cap K = \emptyset ,$$ $$  \vert \eta^{(k)} - \mu \vert  (\Omega) < \frac{1}{4^{k}} , \text{ }  \vert \gamma^{(k)} - v \vert  (\Omega) < \frac{1}{4^{k}} ,$$ 
	where  $  \gamma^{(k)} =\lambda_{k}  \tilde{T}_{\alpha,w}^{* n_{k}} (\eta^{(k)}) $ 
	for each  $ k \in \mathbb{N} ,$ that is  
	$$   \gamma^{(k)} (E) = \int_{E} \lambda_{k} \prod ^{n_{k}-1}_{j=0}( w \circ \alpha^{j-n_{k}} ) d\eta^{(k)} \circ \alpha^{-n_{k}} $$ 
	for every measurable subset  $ E $  of  $ \Omega $  and for each  $ k\in \mathbb{N}  .$  Then we also have that  
	
	$$ \frac{1}{\lambda_{k}}   \tilde{S}_{\alpha,w}^{*n_{k}} ( \gamma^{(k)}) = \eta^{(k)}	
	,$$
	that is  
	$$  \int_{E}  \frac{1}{\lambda_{k}}  \prod ^{n_{k}}_{j=1}( w \circ \alpha^{n_{k}-j} )^{-1}  d\gamma^{(k)} \circ \alpha^{n_{k}  } = \eta^{(k)} (E)  $$ 
	for every measurable subset $E$ of  $ \Omega $  and each  $ k \in \mathbb{N}.  $ Thus we have that
	$$ \vert \lambda_{k} \vert  \vert \tilde{T}_{\alpha,w}^{* n_{k}} (\eta^{(k)}) \vert (K^{c})=
	\vert  \lambda_{k}   \tilde{T}_{\alpha,w}^{* n_{k}} (\eta^{(k)})  \vert  (K^{c}) =
	\vert  \lambda_{k}   \tilde{T}_{\alpha,w}^{* n_{k}} (\eta^{(k)}) - v  \vert  (K^{c}) $$
	$$
	\leq \vert  \lambda_{k}   \tilde{T}_{\alpha,w}^{* n_{k}} (\eta^{(k)}) - v  \vert  (\Omega) <   \frac{1}{4^{k}}
	$$
	and, similarly, we have that
	$$\frac{1}{\vert   \lambda_{k}   \vert}   
	\vert   \tilde{S}_{\alpha,W}^{* n_{k}} (\gamma ^{(k)})  \vert  ( K^{c} )  =
	\vert   \eta^{(k)}   \vert  ( K^{c} ) < \frac{1}{4^{k}}$$ 
	for each $ k \in \mathbb{N}, $ so, for all  $  k \in \mathbb{N} $ we obtain that
	$$
	\vert \tilde{T}_{\alpha,w}^{* n_{k}} (\eta^{(k)})  \vert  ( K^{c} )  
	<  \frac {1 }{  \vert   \lambda_{k}   \vert 4^{k} } , \text{ }	\vert \tilde{S}_{\alpha,W}^{* n_{k}} (\gamma^{(k)})  \vert  ( K^{c} )  
	<  \frac {\vert   \lambda_{k}   \vert }{   4^{k} }  .$$
	By exactly the same arguments as in the proof of \cite[Proposition 3.1]{Grada}, we can deduce that
	$$ \vert  \eta_{}^{(k)} \vert  (A_{k}) <  \frac{1}{2^{k-1}}, \text{ }  \vert  \gamma_{}^{(k)} \vert  (B_{k}) <  \frac{1}{2^{k}-1} ,$$  
	where for each  $ k \in \mathbb{N} $ we put
	$$A_{k} = \lbrace t \in K \text{ }\vert     \prod ^{n_{k}-1}_{i=0}( w \circ \alpha^{i}) (t)  > \frac {1 }{  \vert   \lambda_{k}   \vert 4^{k} }          \rbrace , \text{ } B_{k} = \lbrace t \in K \text{ }\vert     \prod ^{n_{k}}_{i=1}( w \circ \alpha^{-i})^{-1} (t)  >  \frac {\vert   \lambda_{k}   \vert }{   4^{k} }          \rbrace   .$$
	This gives for all  $  k \in \mathbb{N} $ that
	$$\vert \mu \vert (A_{k}) < \frac{1}{2^{k-1}}  + \frac{1}{4^{k}}, \text{ } \vert v \vert (B_{k}) < \frac{1}{2^{k-1}}  + \frac{1}{4^{k}} .$$ 
	Moreover, we have for all $  k \in \mathbb{N} $ that 
	$$ 
	\left( \sup_{t \in K \cap A^{c}_{k}}\prod ^{n_{k}-1}_{i=0}\left( w \circ \alpha^{i} \right) ( t ) \right) \cdot      
	\left( \sup_{t \in K \cap B^{c}_{k}}
	\prod ^{n_{k}}_{i=1}\left( w \circ \alpha^{-1}\right)^{-1} ( t ) \right) < \frac{1}{16^{k}} ,$$
	which proves the implication  $i) \Rightarrow ii).  $ 
	
	Next ve prove  $ ii) \Rightarrow i) .$   Given two non-empty open subsets $\mathcal{O}_{1} $ and $\mathcal{O}_{2} $ of $M(\Omega)  ,$ since $  \mathcal{O}_{1} \setminus \lbrace 0 \rbrace  $ and 
	$ \mathcal{O}_{2} \setminus \lbrace 0 \rbrace $ are also open and non-empty,  we can find some  $ \mu \in \mathcal{O}_{1} \setminus \lbrace 0 \rbrace  $ and some
	$v \in \mathcal{O}_{2} \setminus \lbrace 0 \rbrace $ such that $ \vert \mu \vert (K^{c}) =  \vert v \vert (K^{c})= 0 $ for some compact subset  $ K $  of $ \Omega .$  Choose the strictly increasing sequence  $ \lbrace n_{k} \rbrace_{k} \in \mathbb{N},  $  and the sequences  $ \lbrace A_{k} \rbrace_{k}, \lbrace B_{k} \rbrace_{k} $ of
	Borel subsets of  $ K $  satisfying the assumptions of $ii)$ with respect to  $\mu,v  $  and $K.$  For each  $ k \in \mathbb{N}  ,$  let  $ \tilde{\mu_{k}},  \tilde{v_{k}} $  be the measures in  $M(\Omega)  $  given by  $ \tilde{\mu_{k}} (E)  =  \mu ( E \cap A_{k}^{c})  $  and  $  \tilde{v_{k}} (E)  =  v ( E \cap B_{k}^{c})  $ 
	for every measurable subset $E$ of  $ \Omega .$ By the same arguments as in the proof of \cite[Proposition 3.1]{Grada}, we can deduce that for all $ k \in \mathbb{N} $ it holds that
	$$ \parallel   \tilde{T}_{\alpha,w}^{*n_{k}} ( \tilde{\mu_{k}}  ) \parallel \leq ( \sup_{t \in A^{c}_{k} \cap K}\prod ^{n_{k}-1}_{j=0}\left( w \circ \alpha^{j} \right) ( t ) ) \parallel \mu \parallel , \text{   (1)   }  $$
	
	$$ \parallel   \tilde{S}_{\alpha,w}^{*n_{k}} ( \tilde{v_{k}}  ) \parallel \leq ( \sup_{t \in B^{c}_{k} \cap K}\prod ^{n_{k}}_{j=1}\left( w \circ \alpha^{-j} \right)^{-1} ( t ) ) \parallel v \parallel  . \text{   (2)   }$$
	
	Since  
	$$ \lim_{k \to \infty } ( \mu - \tilde{\mu}_{k} ) =  \lim_{k \to \infty }  (v - \tilde{v}_{k} )= 0 ,$$  
	we may without loss of generality assume that  $ \tilde{\mu}_{k} \in \mathcal{O}_{1} \setminus \lbrace 0 \rbrace $  and  $ \tilde{v}_{k} \in \mathcal{O}_{2} \setminus \lbrace 0 \rbrace  $  for all  $ k \in \mathbb{N} .$ Then, $ \tilde{T}_{\alpha,w}^{*}(\tilde{\mu_{k}}) \neq 0 $ and $ \tilde{S}_{\alpha,w}^{*}(\tilde{v_{k}}) \neq 0 $ for all $  k \in \mathbb{N} $ since $ \tilde{T}_{\alpha,w}^{*} $ and $ \tilde{S}_{\alpha,w}^{*}$ are invertible. For each  $ k \in \mathbb{N} ,$  set
	$$ \eta_{k} = \tilde{\mu_{k}} + \frac{ \parallel \tilde{T}_{\alpha,w}^{* n_{k}} (\tilde{\mu_{k}})  \parallel^{\frac{1}{2}}} 
	{ \parallel \tilde{S}_{\alpha,w}^{* n_{k}} (\tilde{v_{k}})  \parallel^{\frac{1}{2}} } \tilde{S}_{\alpha,w}^{* n_{k}} (\tilde{v_{k}})  .$$
	
	By combining $(1)$ and $(2)$ together with the arguments from the proof of the implication  $ ii) \Rightarrow i) $  in Proposition \ref{proposition-solid}, we can deduce from the assumptions in $ ii)$ that  $ \eta_{k} \rightarrow \mu $  and
	$$
	\frac
	{ \parallel \tilde{S}_{\alpha,w}^{* n_{k}} (\tilde{v_{k}})  \parallel^{\frac{1}{2}}
	} 
	{	\parallel \tilde{T}_{\alpha,w}^{* n_{k}} (\tilde{\mu_{k}})  \parallel^{\frac{1}{2}}  } 
	\tilde{T}_{\alpha,w}^{* n_{k}} (\eta^{(k)}) \rightarrow v $$
	in $M(\Omega)$ as $k \rightarrow \infty,$
	so  $ \tilde{T}_{\alpha,w}^{*} $ is topologically semi-transitive on  $M(\Omega)  .$  
\end{proof}
\begin{proposition}\label{measures-cesaro}
	The following statements are equivalent.\\
	i)  $\tilde{ T }_{\alpha,w}^{*}  $ is topologically Ces\'{a}ro hyper-transitive on  $ M(\Omega) .$ \\
	ii) For every compact subset  $ K $ of $ \Omega $ and any two measures  $\mu, v  $ in  $ M(\Omega)  $  with  $ \vert \mu \vert  (K^{c}) = \vert v \vert  (K^{c}) =0 $ there exist a strictly increasing sequence  $ \lbrace n_{k} \rbrace_{k}  \subseteq \mathbb{N} $  and sequences  $ \lbrace A_{k} \rbrace_{k} , \lbrace B_{k} \rbrace_{k} $ of Borel subsets of $K$ such that
	$ \alpha^{n_{k}}(K) \cap K = \varnothing $  for all  $ k \in \mathbb{N}  $  and
	$$\lim_{k \to \infty } \vert \mu \vert (A_{k}) = \lim_{k \to \infty } \vert v \vert (B_{k})=0,  $$
	
	$$ \lim_{k \to \infty } 
		\left( \sup_{t \in A^{c}_{k} \cap K} n_{k}^{-1}\prod ^{n_{k}-1}_{j=0}\left( w \circ \alpha^{j}\right) ( t ) \right) =  \lim_{k \to \infty }     
	\left( \sup_{t \in B^{c}_{k} \cap K}
	 n_k \prod ^{n_{k}}_{j=1}\left( w \circ \alpha^{-j}\right)^{-1} ( t ) \right)  =0 .$$ 
	 Moreover, if $ \alpha $  is not aperiodic, then we only have $ ii) \Rightarrow i) .$
\end{proposition}
\begin{proof}
	For the proof of the implication $  i) \Rightarrow  ii)  $ we can proceed in the similar way as in the proof of $  i) \Rightarrow  ii)  $ in Proposition \ref{radon} by letting $ \lambda_{k} = n_{k}^{-1} ,$ whereas for the proof of the opposite implication, we can proceed as in the proof of  \cite[Proposition 3.1]{Grada} by considering $ n^{-1}  \tilde{ T }_{\alpha,w}^{*n} $ and $ n  \tilde{ S }_{\alpha,w}^{*n} $  instead of $  \tilde{ T }_{\alpha,w}^{*n} $ and $   \tilde{ S }_{\alpha,w}^{*n} ,$  respectively. 
\end{proof}	
\begin{example}\label{measures}
	Let $ \Omega = \mathbb{R}, \text{ } \alpha : \mathbb{R} \rightarrow \mathbb{R}   $  be given by  $\alpha (t) = t+1  $ for all  $ t \in \mathbb{R}$ and $ w $  be a continuous bounded positive function on  $ \mathbb{R} .$  
	If there exist some  $ M, \delta , K_{1} , K_{2} > 0 $  such that  $ 1 < M-\delta \leq w(t) \leq M $  for all  $ t \leq -K_{1}  $  and  $ w(t) = 1  $  for all  $t \geq K_{2}  ,$  then the conditions of Proposition \ref{measures-cesaro} are satisfied, however, by \cite[Proposition 3.1]{Grada} it follows that  $\tilde{ T }_{\alpha,w}^{*}  $ is not topologically transitive on  $ M(\Omega) .$ On the other hand, if  $ \alpha (t) = t-1  $ for all $ t \in \mathbb{R} , \frac{1}{M} \leq w(t) \leq \frac{1}{M- \delta}  $ for all  $t \leq -K_{1} $ and $ 1 =w(t) $ for all $ t\geq K_{2},$ then the conditions of Proposition \ref{radon} are satisfied, however,  $\tilde{ T }_{\alpha,w}^{*}  $ is neither topologically Ces\'{a}ro hyper-transitive nor topologically transitive on  $ M(\Omega) ,$ which follows from Proposition \ref{measures-cesaro} and \cite[Proposition 3.1]{Grada}. 
\end{example}
\section{Porosity and weighted composition operators}

We recall first the following definition. 
\begin{definition}
	Let $0<\lambda<1$. A subset $E$ of a metric space $X$ is called \emph{$\lambda$-porous} at $x\in E$ if for each $\delta>0$ there is an element $y\in B(x;\delta)\setminus\{x\}$ such that
	$$B(y;\lambda\,d(x,y))\cap E=\varnothing.$$
	$E$ is called \emph{$\lambda$-porous} if it is $\lambda$-porous at every element of $E$. Also, $E$ is called \emph{$\sigma$-$\lambda$-porous} if it is a countable union of $\lambda$-porous subsets of $X$.
\end{definition}
The following lemma will play a key role in the proof of the next proposition in this paper, see also \cite[Lemma 2]{bay}.
\begin{lemma}\label{lem1}
	Let $\mathcal F$ be a non-empty family of non-empty closed subsets of a complete metric space $X$ such that for each $F\in\mathcal F$ and each $x\in X$ and $r>0$ with $B(x;r)\cap F\neq \varnothing$, there exists an element $J\in\mathcal F$ such that 
	$$\varnothing\neq J\cap B(x;r)\subseteq F\cap B(x;r)$$
	and $F\cap B(x;r)$ is not $\lambda$-porous at all elements of $J\cap B(x;r)$. Then, every set in $\mathcal F$ is not $\sigma$-$\lambda$-porous.
\end{lemma}
In the sequel, we shall consider an arbitrary number $0<\lambda\leq \frac{1}{2}$ and for the simplicity, we shall just  write $\sigma$-porous instead of $\sigma$-$\lambda$-porous for a general $ \lambda $ with $0<\lambda\leq \frac{1}{2} .$\\
The proof of the next theorem is motivated by the proof of \cite[Theorem 1]{bay} and \cite[Theorem 2.3]{taiwanese}. 

\begin{theorem}\label{porous}  
	For each $ g \in C_0(\mathbb{R}) ,$ the set $$ \Gamma_g := \{ f \in C_0(\mathbb{R}) \mid |f(m)| \geq |g(m)| \text{ for all  } m \in \mathbb{Z} \} $$ is not $ \sigma - $ porous in $ C_0(\mathbb{R}). $ 
\end{theorem}

\begin{proof} 
	Let $ 0 \leq \lambda \leq \frac{1}{2}  $ and fix some $ 0 < \beta < \lambda .$ 	Put
	$$\mathcal F:=\big\{\Gamma_{g}:\,\, g\in C_0(\mathbb{R}) \big\}.$$
	We will show that the collection $\mathcal F$ satisfies the conditions of Lemma \ref{lem1}.
	Let $g \in C_0(\mathbb{R}) .$ Without loss of generality, we can assume that g is a nonnegative function. Obviously, $ \Gamma_g $ is closed, non-empty subset of $C_0(\mathbb{R}).$  Let $ f \in C_0(\mathbb{R})$ and $ \tilde{r}>0 $ be such that $B(f, \tilde{r}) \cap \Gamma_g \neq \emptyset.$ Our aim in the rest of this proof, as in the proof of \cite[Theorem 1]{bay} and \cite[Theorem 2.3]{taiwanese}, will be to find a nonnegative function $ h \in C_0(\mathbb{R}) $ such that $\emptyset \neq B(f, \tilde{r}) \cap \Gamma_h \subseteq B(f, \tilde{r}) \cap \Gamma_g.$ Let $k \in B(f, \tilde{r}) \cap \Gamma_g $ and $  r \in (0, \tilde{r} - \|k - f\|_{\infty}). $
	Since, $k, f, \beta^{-1}g \in C_0(\mathbb{R}),$ there exists some $ N \in \mathbb{N} $ such that   $$|k(t)|,|f(t)|,\beta^{-1} g(t) < \frac{r}{6}   $$ for all $t \in \mathbb{R} $ with $ |t| \geq N .$ We let $ \delta \in \left(0, \frac{r}{100}\right) $ and we define the function 
	$$ h(t) =
	\begin{cases} 
		g(t) + \delta  \text{ for } t \in [-N, N], \\
		\beta^{-1} g(t)  \text{ for } t \in (-\infty, -N-1] \cup [N+1, \infty), \\
		\beta^{-1} g(t-N-1) + (t + N + 1)\left(\delta + g(-N) - \beta^{-1} g(-N-1)\right) , \text{ for } t \in (-N-1, -N), \\
		g(N) + \delta + (t - N)\left(\beta^{-1} g(N+1) - \delta - g(N)\right)  \text{ for } t \in (N, N+1).
	\end{cases}
	$$
	Then $  h \in \Gamma_g .$ To simplify notation, we set $ M = 2N$  and for each  $j \in \mathbb{Z} \cap [-N, N]$, we put $x_j = j - N$  for  $j \in \{0, 1, \dots, M\}.$ We let then $ \eta $ be the function from $ \{x_0, x_1, \dots, x_M\} $ into $ \mathbb{R} $ given by $\eta(x_i) =
	\begin{cases}
		\frac{k(x_i)}{|k(x_i)|} & \text{if } k(x_i) \neq 0, \ i \in \{0, 1, \dots, M\}, \\
		1 & \text{if } k(x_i) = 0, \ i \in \{0, 1, \dots, M\}.
	\end{cases}
	$
	Further,  we let $\tilde{\eta}$ be the piecewise linear function on $ [-N, N] $ connecting the points $ \left(x, \eta(x_i)\right) $ where $ i \in \{0, \dots, M\}.$ More precisely, on each segment $ [x_{i-1}, x_i] $ with $ i \in \{1, \dots, M\}, $ the function $ \tilde{\eta}$ is given by  $
	\tilde{\eta}(t) = \eta(x_{i-1}) + \frac{t - x_{i-1}}{x_i - x_{i-1}} \left(\eta(x_i) - \eta(x_{i-1})\right) $ for $ t \in [x_{i-1}, x_i].$\\
	Finally, we construct the function $ \mathcal{E} : \mathbb{R} \rightarrow \mathbb{C}
	$ by 
	\\$\mathcal{E}(t) =
	\begin{cases} 
		k(t) + \delta \tilde{\eta}(t)  \text{ for } t \in [-N, N], \\
		h(t)  \text{ for } t \in (-\infty, -N-1] \cup [N+1, \infty), \\
		h(-N-1) + (t+N+1)\left(k(-N) + \delta \tilde{\eta}(-N) - h(-N-1)\right)  \text{ for } t \in (-N-1, -N), \\
		k(N) + \delta \tilde{\eta}(N) + (t-N)\left(h(N+1) - k(N) - \delta \tilde{\eta}(N)\right)  \text{ for } t \in (N, N+1).
	\end{cases}
	.$
	Then $ \mathcal{E} \in C_0(\mathbb{R}) $ since $ h \in C_0(\mathbb{R}) ,$ and $  k, \tilde{\eta}$  are continuous.  We notice that by the triangle inequality, for all $ t \in (N, N+1) $ we have $$ |\mathcal{E}(t)| = |(N+1-t)(k(N) + \delta \tilde{\eta}(N)) + (t-N)h(N+1)| $$ $$
	\leq (N+1-t)(|k(N)| + \delta) + (t-N)\beta^{-1}g(N+1) \leq (N+1-t)\left(\frac{r}{6} + \frac{r}{100}\right) + (t-N)\frac{r}{6} $$ $$ \leq \frac{r}{6} + \frac{r}{100} < \frac{r}{3} .$$
	Similarly, we have for all $t \in (-N-1, -N] $ that $ |\mathcal{E}(t)| < \frac{r}{3}  .$ Therefore, $$ |\mathcal{E}(t) - f(t)| \leq |\mathcal{E}(t)| + |f(t)| \leq \frac{r}{3} + \frac{r}{6} = \frac{r}{2} < \tilde r. $$ for all $ t \in (-N-1, -N) \cup (N, N+1) .$ Let now $ t\in [-N, N].$ Then there exists some $ i \in \{1, \dots, M\} $ such that $ t \in [x_{i-1}, x_i],$ hence, we obtain that $$|k(t) - \mathcal{E}(t)| = \delta |\tilde{\eta}(t)| = \delta \left|\eta(x_{i-1}) + \frac{t - x_{i-1}}{x_i - x_{i-1}}(\eta(x_i) - \eta(x_{i-1}))\right| $$ $$
	\leq \delta \left[\left(1 - \frac{t - x_{i-1}}{x_i - x_{i-1}}\right)|\eta(x_{i-1})| + \frac{t - x_{i-1}}{x_i - x_{i-1}}|\eta(x_i)|\right] \leq \delta.$$ This holds for all $t \in [-N, N].$ Therefore, for all $ t \in [-N, N]$ we obtain 
	$$|\mathcal{E}(t) - f(t)| \leq |\mathcal{E}(t) - k(t)| + \|k - f\|_{\infty} < \delta + \|k - f\|_{\infty} < r + \|k - f\|_{\infty} < \tilde{r} $$ 
	Finally, for all $t \in (-\infty, -N-1] \cup [N+1, \infty),$ we get that $$
	|\mathcal{E}(t) - f(t)| \leq |\mathcal{E}(t)| + |f(t)| = \beta^{-1}g(t) + |f(t)| < \frac{r}{3} < \tilde{r}. 
	$$ Hence, $ \|\mathcal{E} - f\|_{\infty} < \tilde{r}. $ Moreover, for all $ m \in \mathbb{Z} \cap [-N, N], $ we have that $$ |\mathcal{E}(m)| = |k(m) + \delta \eta(m)| = |k(m)| + \delta,
	$$ since $ \tilde{\eta}(m) = \eta(m) 	$ for all $ m \in \mathbb{Z} \cap [-N, N].$ Hence, we deduce that $$ |\mathcal{E}(m)| = |k(m)| + \delta \geq g(m) + \delta = h(m), $$ and for each $ m \in \mathbb{Z} \cap (-\infty, -N-1] \cup [N+1, \infty), 
	|\mathcal{E}(m)| = h(m), $ so $ \mathcal{E} \in \Gamma_h. $ Thus, $ \emptyset \neq B(f, \tilde{r}) \cap \Gamma_h \subseteq B(f, \tilde{r}) \cap \Gamma_g. $
	
	As in the proof of \cite[Theorem 1]{bay} and \cite[Theorem 2.3]{taiwanese}, we let then $ u \in B(f, \tilde{r}) \cap \Gamma_h, r^{'} = \min \left\{\delta, \lambda(\tilde{r} - \|f - u\|_{\infty})\right\}$ and we pick some $ v \in B(u, r^{'}).$ We let $ \Theta $ be the function from 
	$$ \mathbb{Z} \cap \left( (-\infty, -N-1] \cup [N+1, \infty) \right) $$ 
	into $ \mathbb{C} $ given by 
	$ \Theta(m) =
	\begin{cases}
		\frac{v(m)}{|v(m)|} & \text{if } v(m) \neq 0, \\
		1 & \text{if } v(m) = 0,
	\end{cases}
	$ and we let $ \tilde{\Theta}$ be the piecewise linear function from $ (-\infty, -N-1] \cup [N+1, \infty) $ into $ \mathbb{C} $ connecting the points $ (m, \Theta(m)) $ where $ m \in \mathbb{Z} \cap \left( (-\infty, -N-1] \cup [N+1, \infty) \right)  .$ More precisely, on each segment $[m-1, m] $ with $ m \in \mathbb{Z} \cap ((-\infty, -N-1] \cup [N+2, \infty)), $ the function $ \tilde{\Theta} $ is given as $$
	\tilde{\Theta}(t) = \Theta(m-1) + (t+1-m)(\Theta(m) - \Theta(m-1)) $$ for $ t \in [m-1, m]. $ By the triangle inequality, it follows that $|\tilde{\Theta}(t)| \leq 1 \text{ for all } t \in (-\infty, -N-1] \cup [N+1, \infty).$\\
	Finally, we construct the function $ \gamma: \mathbb{R} \rightarrow \mathbb{C}
	$ given by \\ $ \gamma(t) =
	\begin{cases}
		v(t)  \text{ if } t \in [-N, N], \\
		v(t) + \beta |u(t) - v(t)|\tilde{\Theta}(t)  \text{ for } t \in (-\infty, -N-1] \cup [N+1, \infty), \\
		v(t) + (t-N)|\beta | u(N+1) - v(N+1)|\tilde{\Theta}(N+1)  \text{ for } t \in (N, N+1), \\
		v(t) - (t+N)\beta| u(-N-1) - v(-N-1)|\tilde{\Theta}(-N-1)  \text{ for } t \in (-N-1, -N).
	\end{cases}
	$
	Since $\tilde{\Theta}$ is continuous,  $ |\tilde{\Theta}(t)| \leq 1 $  for all $ t \in (-\infty, -N-1] \cup [N+1, \infty) $ and $ u, v \in C_0(\mathbb{R})
	,$ it is not hard to see that $ \gamma \in C_0(\mathbb{R}) .$ Also, from the construction of the function $ \gamma $ it follows that $ \|\gamma - v\|_{\infty} \leq \beta \|u - v\|_{\infty} \leq \lambda \|u - v\|_{\infty} .$ For $ m \in \mathbb{Z} \cap \left( (-\infty, -N-1] \cup [N+1, \infty) \right) $ we have\\
	$ \gamma(m) =
	\begin{cases}
		v(m) \left(1 + \frac{\beta | u(m) - v(m)}{|v(m)|}\right) & \text{if } v(m) \neq 0, \\
		\beta |u(m)| & \text{if } v(m) = 0,
	\end{cases}
	$ because $ \tilde{\Theta}(m) = \Theta(m)$ for all $ m \in \mathbb{Z} \cap \left( (-\infty, -N-1] \cup [N+1, \infty) \right).$ Thus, since $ |u(m)| \geq h(m) $ for all $ m \in \mathbb{Z} ,$ we get $$ |\gamma(m)| = |v(m)| + \beta |u(m) - v(m)| \geq \beta |u(m)| \geq \beta h(m) = g(m) $$ for all $ m \in \mathbb{Z} \cap \left( (-\infty, -N-1] \cup [N+1, \infty) \right) .$  \\
	Moreover, since $ v \in B(u, r^{'}), $ we have $ \|u - v\|_{\infty} \leq r^{'} \leq \delta, $ hence $ \\
	|\gamma(m)| = |v(m)| \geq |u(m)| - \delta \geq h(m) - \delta = g(m) $ for all $ m \in \mathbb{Z} \cap [-N, N] .$ Therefore, $B(v, \lambda \|u - v\|_{\infty}) \cap B(f, \tilde{r}) \cap \Gamma_g \neq \emptyset. $
\end{proof}
We obtain the following corollary of Theorem \ref{porous}.

\begin{corollary}\label{porosity}
	Consider the weighted composition operator $\tilde{T}_{\alpha, w} $ on $ C_0(\mathbb{R}) $ given by $ \tilde{T}_{\alpha, w}(f) = w \cdot (f \circ \alpha), $ where $ 0 < w, w^{-1} \in C_b(\mathbb{R}), $ and $ \alpha $ is a homeomorphism of $ \mathbb{R}. \\
	$ If $ \lim_{n \to \infty} \prod_{k=1}^{n} (w \circ \alpha^{-k})^{-1}(n) = 0, $ then the set $$ 
	\left\{f \in C_0(\mathbb{R}) : \|\tilde{T}_{\alpha, w}^n(f)\|_{\infty} \geq 1 \text{ for all } n \in \mathbb{N}\right\} \\
	$$ is not $ \sigma\text{-porous. In particular, the set of non-hypercyclic vectors for the operator } \tilde{T}_{\alpha, w} $ is not $ \sigma\text{-porous in } C_0(\mathbb{R}).$
\end{corollary}	

\begin{proof}
	We let $ g $ be the piecewise linear function on $\mathbb{R}^{+} $ connecting the points $ \left(n, \prod_{k=1}^{n} (w \circ \alpha^{-k})^{-1}(n)\right)$ where $ n \in \mathbb{N}. $ More precisely, for each $ n \in \mathbb{N}, $ we let 
	$$g(t) = \prod_{k=1}^{n} (w \circ \alpha^{-k})^{-1}(n)$$ $$ + (t-n) \left(\prod_{k=1}^{n+1} (w \circ \alpha^{-k})^{-1}(n+1) - \prod_{k=1}^{n} (w \circ \alpha^{-k})^{-1}(n)\right) $$ for $ t \in [n, n+1], $ whereas for $ t \in [0, 1], $ we put $ g(t) = t (w \circ \alpha^{-k})^{-1}(1). $ Moreover, we let $ g(t) = 0 $ for $ t \leq 0.$ It is easily seen that $ g \geq 0 $ and $ g \in C_0(\mathbb{R}) $ since $ \lim_{n \to \infty} \prod_{k=1}^{n} (w \circ \alpha^{-k})^{-1}(n) = 0 $ by the assumption. By Theorem \ref{porous} , the set $\Gamma_g $ is not $ \sigma\text{-porous.} $ Now, for each $ f \in \Gamma_g $ and $ n \in \mathbb{N}, $ we have that $$  \|\tilde{T}_{\alpha, w}^n(f)\|_{\infty} = \left\|\prod_{k=1}^{n} (w \circ \alpha^{n-k}) (f \circ \alpha^n)\right\|_{\infty} = \left\|\prod_{k=1}^{n} (w \circ \alpha^{-k}) f\right\|_{\infty}$$ $$ \geq \prod_{k=1}^{n} (w \circ \alpha^{-k})(n) |f(n)| \geq \prod_{k=1}^{n} (w \circ \alpha^{-k})(n) g(n) = 1. $$ \end{proof}

\bibliographystyle{amsplain}

\end{document}